\input harvmac.tex 
\newdimen\tableauside\tableauside=1.0ex
\newdimen\tableaurule\tableaurule=0.4pt
\newdimen\tableaustep
\def\phantomhrule#1{\hbox{\vbox to0pt{\hrule height\tableaurule width#1\vss}}}
\def\phantomvrule#1{\vbox{\hbox to0pt{\vrule width\tableaurule height#1\hss}}}
\def\sqr{\vbox{%
  \phantomhrule\tableaustep
  \hbox{\phantomvrule\tableaustep\kern\tableaustep\phantomvrule\tableaustep}%
  \hbox{\vbox{\phantomhrule\tableauside}\kern-\tableaurule}}}
\def\squares#1{\hbox{\count0=#1\noindent\loop\sqr
  \advance\count0 by-1 \ifnum\count0>0\repeat}}
\def\tableau#1{\vcenter{\offinterlineskip
  \tableaustep=\tableauside\advance\tableaustep by-\tableaurule
  \kern\normallineskip\hbox
    {\kern\normallineskip\vbox
      {\gettableau#1 0 }%
     \kern\normallineskip\kern\tableaurule}%
  \kern\normallineskip\kern\tableaurule}}
\def\gettableau#1 {\ifnum#1=0\let\next=\null\else
  \squares{#1}\let\next=\gettableau\fi\next}

\tableauside=1.0ex
\tableaurule=0.4pt

\lref\cs{E. Witten, ``Quantum field theory and the Jones polynomial," 
Commun. Math. Phys. {\bf 121} (1989) 351.}
\lref\gv{R. Gopakumar and C. Vafa, ``On the gauge theory/geometry 
correspondence," hep-th/9811131, Adv. Theor. Math. Phys. {\bf 3} (1999) 1415.}
\lref\ov{H. Ooguri and C. Vafa, ``Knot invariants and topological 
strings," hep-th/9912123, Nucl. Phys. {\bf B 577} (2000) 419.} 
\lref\gvmone{R. Gopakumar and C. Vafa, ``M-theory and topological 
strings, I," hep-th/9809187. }
\lref\gvm{R. Gopakumar and C. Vafa, ``M-theory and topological 
strings, II," hep-th/9812127. }
\lref\kkv{S. Katz, A. Klemm and C. Vafa, ``M-theory, topological strings, and 
spinning black-holes," hep-th/9910181, 
Adv. Theor. Math. Phys. {\bf 3} (1999) 1445.}
\lref\fh{W. Fulton and J. Harris, {\it Representation theory. A first course},
Springer-Verlag, 1991.}
\lref\jonesann{V.F.R. Jones, ``Hecke algebras representations of braid 
groups and link polynomials," Ann. of Math. {\bf 126} (1987) 335. }
\lref\guada{E. Guadagnini, ``The universal link polynomial," 
Int. J. Mod. Phys. {\bf A 7} (1992) 877;  {\it The link invariants of the 
Chern-Simons field theory,} Walter de Gruyter, 1993.}
\lref\lick{W.B.R. Lickorish, {\it An introduction to knot theory}, 
Springer-Verlag, 1998.}
\lref\awrep{M. Wadati, T. Deguchi and Y. Akutsu, ``Exactly solvable models 
and knot theory," Phys. Rep. {\bf 180} (1989) 247.}
\lref\ofer{O. Aharony, S. Gubser, J. Maldacena, H. Ooguri and Y. Oz, 
``Large $N$ field theories, string theory and gravity," hep-th/9905111, 
Phys. Rep. {\bf 323} (2000) 183.}
\lref\homfly{P. Freyd, D. Yetter, J. Hoste, W.B.R. Lickorish, K. Millet and 
A. Ocneanu, ``A new polynomial invariant of knots and links," Bull. Amer. 
Math. Soc. {\bf 12} (1985) 239.}
\lref\lpp{J.M.F. Labastida and E. P\'erez, ``Gauge-invariant 
operators for singular knots in
Chern-Simons gauge theory," hep-th/9712139, Nucl.\ Phys.\ {\bf
B 527} (1998) 499.}
\lref\witop{E. Witten, ``Chern-Simons gauge theory as 
a string theory,'' hep-th/9207094, in {\it The Floer memorial volume}, 
H. Hofer, C.H. Taubes, A. Weinstein and E. Zehner, eds., 
Birkh\"auser 1995, p. 637.}  
\lref\douglas{M.R. Douglas, ``Chern-Simons-Witten theory as a topological 
Fermi liquid,'' hep-th/9403119.}
\lref\lm{J.M.F. Labastida and M. Mari\~no, ``Polynomial invariants 
for torus knots and topological strings,''  hep-th/0004196,
Commun. Math. Phys. {\bf 217} (2001) 423.}
\lref\alp{M. \'Alvarez, J.M.F. Labastida, and E. P\'erez, 
``Vassiliev invariants for links from Chern-Simons perturbation theory,'' 
hep-th/9607030, Nucl. Phys. {\bf B 488} (1997) 677.}
\lref\lickm{W.B.R. Lickorish and K.C. Millett, ``A polynomial invariant of 
oriented links,'' Topology {\bf 26} (1987) 107.}
\lref\cp{V. Chari and A. Presley, {\it A guide to quantum groups}, 
Cambridge University Press, 1994.}
\lref\tur{V.G. Turaev, ``The Yang-Baxter equation and invariants 
of links,'' Inv. Math. {\bf 92} (1988) 527.}
\lref\rosso{M. Rosso, ``Groupes quantiques et mod\`eles \`a vertex 
de V. Jones en th\'eorie des noeuds,'' C.R. Acad. Sci. Paris {\bf 307} 
(1988) 207.} 
\lref\rj{M. Rosso and V. Jones, ``On the invariants of torus knots 
derived from quantum groups,'' J. Knot Theory Ramifications {\bf 2} (1993) 
97.}
\lref\lmv{J.M.F. Labastida, M. Mari\~no and C. Vafa, ``Knots, links and 
branes at large $N$,'' hep-th/0010102, JHEP {\bf 0011} (2000) 007.}
\lref\kl{S. Katz and M. Liu, ``Enumerative geometry of stable 
maps with Lagrangian boundary conditions and multiple covers of the disc,'' 
math.AG/0103074.}
\lref\ls{J. Li and Y.S. Song, ``Open string instantons and relative 
stable morphisms,'' hep-th/0103100.}
\lref\vass{V. A. Vassiliev, ``Cohomology of knot spaces", {\it Theory
of singularities and its applications}, Advances in Soviet
Mathematics, vol. 1, Americam Math. Soc., Providence, RI, 1990,
23-69.}
\lref\birrev{J.S. Birman, ``New points of view in knot theory", Bull. Amer. 
Math. Soc. {\bf 28} (1993) 253.}
\lref\barnatan{D. Bar-Natan, ``On the Vassiliev knot invariants", Topology
{\bf 34} (1995) 423.}
\lref\bilin{J.S. Birman and X.S. Lin, ``Knot polynomials and Vassiliev
invariants", Inv. Math. {\bf 111} (1993) 225.}
\lref\natan{D. Bar-Natan, ``Perturbative aspects of Chern-Simons
topological quantum field theory", Ph.D. Thesis, Princeton University,
1991.}
\lref\rama{P. Ramadevi and T. Sarkar, ``On link invariants and 
topological string amplitudes,'' hep-th/0009188, Nucl. Phys. {\bf B 600} 
(2001) 487.}
\lref\kau{L.H. Kauffman, ``An invariant of regular isotopy", 
Trans. Amer. Math. Soc. {\bf 318} (1990) 417.}
\lref\kont{M. Kontsevich, ``Vassiliev knot invariants", Advances in Soviet
Math. {\bf 16}, Part 2 (1993) 137.}
\lref\aires{J.M.F.  Labastida, ``Chern-Simons gauge theory: ten years after",
hep-th/9905057, in {\sl Trends in Theoretical Physics II}, H. Falomir, R. Gamboa,  F.
Schaposnik, eds., American Institute of Physics, New York, 1999, CP 484,
1-41.}
\lref\periwal{V. Periwal, ``Topological closed string 
interpretation of Chern-Simons theory,'' hep-th/9305115, Phys. Rev. Lett. 
{\bf 71} (1993) 1295. }
\lref\sv{S. Sinha and C. Vafa, ``SO and Sp Chern-Simons at large $N$,'' 
hep-th/0012136.}
\lref\thooft{G. 't Hooft, ``A planar diagram theory for strong interactions,'' 
Nucl. Phys. {\bf B 79} (1974) 461.}
\lref\cg{R. Correale and E. Guadagnini, 
``Large-$N$ Chern-Simons field theory,'' Phys. Lett. {\bf B 337} (1994) 80.}
\lref\gmm{E. Guadagnini, M. Martellini and M. Mintchev, ``Wilson 
lines in Chern-Simons theory and link invariants,'' Nucl. Phys. 
{\bf B 330} (1990) 575.}
\lref\singer{S. Axelrod and I.M. Singer, ``Chern-Simons Perturbation
Theory", 1991, hep-th/9110056 and J. Diff. Geom. {\bf 39} (1994) 173.}
\lref\bt{R. Bott and C. Taubes,  ``On the self-linking of knots", Jour.
Math. Phys. {\bf 35} (1994) 5247.}
\lref\gk{E. Getzler and M.M. Kapranov, ``Modular operads,'' dg-ga/9408003, 
Compositio Math. {\bf 110} (1998) 65.}  
\lref\bp{J. Bryan and R. Pandharipande, ``BPS states of curves in 
Calabi-Yau three-folds,'' math.AG/0009025, Geom. Topol. {\bf 5} (2001) 287.} 

\Title{\vbox{\baselineskip12pt
\hbox{US-FT-4/01}
\hbox{RUNHETC-2001-11}
\hbox{math.QA/0104180}
}}
{\vbox{\centerline{A New Point of View in the Theory}
\centerline{ }
\centerline{of Knot and Link Invariants}}
}
\centerline{\bf Jos\'e M. F. Labastida$^{a}$ and Marcos Mari\~no$^{b}$}

\bigskip
\medskip
{\vbox{\centerline{$^{a}$ \sl Departamento de F\'\i sica de Part\'\i culas}
\centerline{\sl Universidade de Santiago de Compostela}
\centerline{\sl E-15706 Santiago de Compostela, Spain}}
\centerline{ \it labasti@fpaxp1.usc.es}

\bigskip
\medskip
{\vbox{\centerline{$^{b}$ \sl New High Energy Theory Center}
 \centerline{\sl Rutgers University}
\centerline{\sl Piscataway, NJ 08855, USA }}
\centerline{ \it marcosm@physics.rutgers.edu }

\bigskip
\bigskip
\noindent
Recent progress in string theory has led to a 
reformulation of quantum-group polynomial invariants for
knots and links into new polynomial invariants whose coefficients can be
understood in topological terms. We describe in detail how to 
construct the new polynomials and we conjecture their general structure. 
This leads to new conjectures on the algebraic structure of the 
quantum-group polynomial invariants. We also describe 
the geometrical meaning of the coefficients in terms
of the enumerative geometry of Riemann surfaces with boundaries 
in a certain Calabi-Yau threefold.

\bigskip
\bigskip
\bigskip
\bigskip

\Date{April, 2001}


\newsec{Introduction}

During the last two decades the theory of knot and link invariants has
experienced important progress. In the eighties a series of new
polynomial invariants were discovered \jonesann\homfly\kau\awrep\ leading to
a unified picture provided by quantum-group polynomial invariants \tur,
and by vacuum expectation values of operators in Chern-Simons gauge
theory \cs. The progress in the nineties was characterized by the
discovery of Vassiliev invariants \vass\ or invariants of finite type. It was
realized soon that these types of invariants were related. It is now known that
the coefficients of the power series expansion of the quantum-group invariants
\bilin\barnatan\ and the coefficients of the perturbative series
expansion of the vacuum expectation values of operators in Chern-Simons
gauge theory \gmm\natan\lpp\ are Vassiliev invariants\foot{For
reviews on these developments see \birrev\ and \aires.}. These connections
inspired important developments in the theory of finite-type invariants
\singer\kont\bt.

Besides the great progress achieved during the last two decades, there are
still many unanswered questions in the theory of knots and links. One of these
questions is about the topological meaning of the integer coefficients of 
the quantum-group polynomial invariants (see, for example, \birrev). 
The discovery of the relation between these polynomial
invariants and Vassiliev invariants did not provide important progress in
this direction. The main goal of this paper is to point out that the
situation has changed dramatically in the last two years. A new point of
view to study knot and link invariants is now available. In this new
approach the integer coefficients of a reformulation of the quantum-group
polynomial invariants carry topological content.

At the heart of this development is the idea that quantum 
gauge theories may have a string theory description. In the case 
of Chern-Simons theory, a first 
step in this direction was taken by Witten in 1992 \witop\foot{This idea 
has been also explored in \periwal\ and \douglas.}. 
The final picture emerged in 1998, when R. Gopakumar and C. Vafa \gv\ 
found a description of 
Chern-Simons theory in terms of a closed, topological string theory. 
In 1999 H. Ooguri and C. Vafa \ov\ showed how to describe Wilson loops 
of knots in Chern-Simons theory by introducing an open-string sector in the 
topological string theory of \gv. They also showed that the 
string description of Chern-Simons Wilson loops involved a reformulation 
of quantum-group invariants in terms of new integer invariants. 
The connection between the quantum-group invariants and the 
string description was spelled out in detail in \lm. The integer 
invariants of \ov\ were further refined in \lmv, 
where the string description was also extended to the case of links. 

In this paper we will present a summary of these developments 
and their mathematical implications for the theory of quantum-group 
invariants of knots and links.  
The new point of view described in this paper first reformulates the
quantum-group polynomial invariants in terms of new polynomials, and then
assigns topological content to their coefficients. The relation establishes a
connection between quantum-group invariants and the geometry of the moduli 
spaces of Riemann surfaces with holes holomorphically embedded in 
a specific Calabi-Yau manifold. The embedding has Lagrangian boundary 
conditions, and the Lagrangian submanifold which specifies them 
is conjectured to be determined by the link \ov. 
This connection has to be considered at the level of conjecture. The integer
invariants are difficult to compute from the topological side and up to date
the conjecture has been fully tested only for the unknot \kl\ls.
Nevertheless, the conjecture also predicts a particular structure for the
reformulated polynomial invariants which has been verified in many cases 
\ov\lm\rama\lmv.

The integer coefficients of the reformulated polynomial invariant of a given 
link can be also interpreted in terms of a generalization of Gromov-Witten 
invariants which involve Riemann surfaces with boundaries. 
These integer invariants turn out to be a resummation of
generalized Gromov-Witten invariants much in the same spirit as
the Gopakumar-Vafa invariants are for the ordinary ones \gvm\kkv. 

The paper is organized as follows. In section 2 we review the polynomial
invariants for knots and links from a quantum group perspective and we
make a slight refinement for the case of links to make them more suitable
for our purposes. In section 3 we reformulate these polynomials in terms of
new ones and we conjecture their general structure. We also 
present some simple consequences of the conjecture 
for the structure of the HOMFLY polynomial of links. In section 4 we
describe the topological content of the
integer coefficients of the new polynomials. Finally, in section 5 we
state our conclusions and comment on future developments.

\newsec{Quantum-group polynomial invariants of knots and links}

The quantum-group invariants that we will be considering 
are multicolored generalizations of the HOMFLY polynomial (with an important 
subtlety in the case of links), and their 
definition is as follows.  

Let ${\cal L}$ be a link of $L$ components ${\cal K}_{\alpha}$, 
$\alpha=1, \cdots, L$. 
The linking number of the components ${\cal K}_{\alpha}$, ${\cal K}_{\beta}$ 
will be denoted by ${\rm lk}({\cal K}_{\alpha}, {\cal K}_{\beta})$. 
We will represent the link ${\cal L}$ by the closure of a braid. The braid
will have ${\cal N}$ strands, and ${\cal N}_{\alpha}$  will denote the number
of strands of the component ${\cal K}_{\alpha}$.   We will also assume that
the strands are ordered in such a way that  the first ${\cal N}_1$ strands
correspond to the component ${\cal K}_1$,  and so on. 


We associate to each component of the link an irreducible 
representation $R_{\alpha}$ of the quantized universal 
enveloping algebras $U_q({\rm sl}(N, {\bf C}))$. These 
representations are labeled by highest weights $\Lambda_{\alpha}$, 
and as usual (see for example \fh) 
we will associate to them a Young diagram with $\ell_{\alpha}$ boxes. 
The corresponding module will be denoted by $V_{\Lambda_{\alpha}}$. 
Therefore, to the $j$-th strand in the braid we will associate an irreducible 
module $V_j$. If the $j$-th strand belongs to the component 
${\cal K}_{\alpha}$, then $V_j =V_{\Lambda_{\alpha}}$. The total 
module associated to the braid is then $V=\otimes_{j=1}^{\cal N} V_j$, which 
is also given by :
\eqn\tenprod{
V=\otimes_{\alpha=1}^L V_{\Lambda_{\alpha}}^{{\cal N}_{\alpha}}.}
 
It is well-known that each solution to the Yang-Baxter equation provides a 
representation of the braid group of ${\cal N}$ 
strands ${\cal B}_{\cal N}$. The solution that we will use is 
the universal ${\cal R}$-matrix of $U_q({\rm sl}(N,{\bf C}))$ \cp:    
\eqn\rmat{
{\cal R}=q^{{1 \over 2}\sum_{i,j} C_{ij}^{-1}H_i \otimes H_j} 
\prod_{\beta}\exp_q[(1-q^{-1}) 
X_{\beta}^+ \otimes X_{\beta}^-]. 
}
In this equation, $C_{ij}$ is the Cartan matrix of $SU(N)$, and the 
$q$-exponential has the form,
\eqn\notat{
\exp_q (x)= \sum_{k=0}^{\infty} q^{ {1 \over 4} k(k+1)}{ x^k \over 
[k]_q!},}
where the $q$-numbers are defined as:
\eqn\qnum{
[n]_q={q^{n\over 2} -q^{-{n\over 2}} \over q^{1 \over 2} -q^{-{1\over 2 }}}.} 
In the expression \rmat\ for ${\cal R}$, the product over $\beta$ denotes 
a product over positive roots of ${\rm sl}(N, {\bf C})$, and the 
$X_{\beta}^{\pm}$ are certain elements in $U_q ({\rm sl}(N, {\bf C}))$ (see 
\cp, chapter 8 for details).
Notice that ${\cal R}$ acts in a natural way on the tensor product 
of two $U_q({\rm sl}(N,{\bf C}))$-modules. 
Finally, we denote $\check{{\cal R}}=P_{12}{\cal R}$, 
where $P_{12}$ is the exchange operator 
between the two factors of the tensor product.

The representation of ${\cal B}_{\cal N}$ on $V$ is defined as follows.
If $\sigma^{\pm 1}_i$ is an elementary braid, then
\eqn\repre{
\pi (\sigma_i^{\pm 1})=I_{V_1} \otimes 
\cdots \otimes  \check{{\cal R}}^{\pm 1} 
\otimes \cdots \otimes I_{V_{\cal N}},}
where $\check{{\cal R}}^{\pm 1}$ acts on $V_i \otimes V_{i+1}$. 
Therefore, every braid word $\xi$ gives an operator on the total 
module $V$ \tenprod\ that we will denote by $\pi (\xi)$. Since 
${\cal R}$ is a solution of the quantum Yang-Baxter equation, the 
above representation of the braid group is well-defined, {\it i.e.}, it 
respects the relations between the generators \tur. 

In order to define an invariant of links, we need 
an enhancement of the ${\cal R}$-matrix \tur. We take \rosso\rj: 
\eqn\enh{
\mu=q^{ \rho^*},}
where $\rho^*$ is the element in the Cartan subalgebra ${\bf h} 
\subset U_q({\bf h})$ which corresponds to the Weyl vector ({\it i.e.}, the 
sum of fundamental weights) under the natural isomorphism ${\bf h} \simeq 
{\bf h}^*$ induced by the Killing form. 

The quantum-group invariant that we will consider is defined as 
follows:
\eqn\qginv{
W_{(R_1, \cdots, R_L)}({\cal L})=q^{d ({\cal L})} 
{\rm Tr}_V \bigl( \mu^{\otimes {\cal N}}\pi(\xi)\bigr),}
where $d({\cal L})$ is given by:
\eqn\corr{
d ({\cal L})=-{1 \over 2} \sum_{\alpha=1}^L w({\cal K}_{\alpha}) 
(\Lambda_{\alpha}, \Lambda_{\alpha} + 2 \rho) + 
{1 \over N}\sum_{\alpha < \beta} {\rm lk}({\cal K}_{\alpha}, 
{\cal K}_{\beta}) \ell_{\alpha} \ell_{\beta}.}
In this expression $w({\cal K}_{\alpha})$ is the writhe of the $\alpha$-th
component  of the link ${\cal L}$, and $\ell_{\alpha}$ and $\ell_{\beta}$ are
the number of boxes of the Young diagrams associated to the irreducible
representations $\Lambda_{\alpha}$ and $\Lambda_{\beta}$.  The first
term in
\corr\ guarantees, by the usual arguments
\tur,  that \qginv\ is an ambient isotopy invariant of the link. The
second term  in \corr\ cancels overall powers of $q^{1/N}$
that appear  after taking the trace in \qginv. The resulting quantum-group
invariant  is in general a {\it rational
function}\foot{Often, following standard usage, we will refer to
this invariant (as well as to the reformulated one below) as polynomial invariant
though, strictly speaking, it is in general a rational function.} of
$q^{\pm {1\over 2}}$ and 
$\lambda^{\pm {1\over 2}}$, where
\eqn\lamb{
\lambda =q^N.}

\vbox{\medskip\noindent{\bf Remarks}:

\noindent$\bullet$ If ${\cal L}$ is the trivial link of $L$ components, with
attached  representations $R_{\alpha}$, $\alpha=1,\cdots,L$, then the
quantum-group invariant is,
\eqn\invtriv{
 W_{(R_1, \cdots, R_L)}({\cal L})= 
\prod_{\alpha=1}^L {\rm dim}_q (R_{\alpha}),}
where ${\rm dim}_q (R_{\alpha})$ is the quantum dimension of $R_{\alpha}$. }

\noindent$\bullet$ When all the components of the link are in the fundamental 
representation, {\it i.e.}, $\tableauside=1.5ex R_{\alpha}=\tableau{1}
\tableauside=1.0ex$,  the above quantum-group invariant is related to the 
HOMFLY polynomial of the link, $P_{\cal L}(q, \lambda)$, in the 
following way:
\eqn\rel{
W_{(\tableau{1}, \cdots,\tableau{1})}({\cal L}) = 
\lambda^{{\rm lk}({\cal L})}
\biggl( {\lambda^{1\over 2} -\lambda^{-{1\over 2}} 
\over q^{1\over 2} -q ^{-{1\over 2}}} \biggr) P_{\cal L}(q, \lambda),}
where 
\eqn\linkn{
{\rm lk}({\cal L}) =\sum_{\alpha<\beta}{\rm lk}({\cal K}_{\alpha}, 
{\cal K}_{\beta})}
is the total linking number of ${\cal L}$. From this relation we observe
that the  quantum-group invariant defined in \qginv\ for a link with 
the fundamental representation attached to all its components is not 
the unnormalized HOMFLY polynomial, but differs from it in an overall factor 
$\lambda^{{\rm lk}({\cal L})}$. Usually, when all the components 
of a link are in the same representation $\Lambda$,
 its quantum-group invariants are defined as in \qginv, 
with the only difference that the 
overall power of $q$ is not $d({\cal L})$ but,
\eqn\usualdef{
-{1 \over 2} (\Lambda, \Lambda + 2 \rho)w({\cal L}),}
where $w({\cal L})$ the total writhe of the link, 
\eqn\writhe{
w({\cal L})=\sum_{\alpha=1}^L w({\cal K}_{\alpha}) + 2\, {\rm lk}({\cal L}).}
The difference between \usualdef\ and \corr\ when
$\tableauside=1.5ex \Lambda=\tableau{1} \tableauside=1.0ex$  gives
precisely the extra factor
$\lambda^{{\rm lk}({\cal L})}$,  which will be crucial for our
considerations.
 
\noindent$\bullet$ In \qginv, we have assumed that none of the representations 
$R_1, \cdots, R_L$ is the trivial one. It will be useful to extend the 
definition to include the trivial representation as follows: let 
$\{ \alpha_1, \cdots, \alpha_s \}$ be a subset of $\{1, \cdots, L\}$, with 
$s>0$. 
The complementary set will be denoted by $\{ \alpha_{s+1}, \cdots, 
\alpha_L \}$. The 
sublink  ${\cal L}_{\alpha_1, 
\cdots, \alpha_s}$ of ${\cal L}$ is obtained by considering only the 
components ${\cal K}_{\alpha_i}$ of ${\cal L}$, with $i=1, \cdots, s$, and 
``deleting'' the rest of the components. If we take 
$R_{\alpha_{s+1}}=\cdots=R_{\alpha_{L}}=\cdot$ to be the 
trivial representation, we define:
\eqn\trivreps{
W_{(R_1, \cdots, R_L)}({\cal L})=W_{(R_{\alpha_1}, \cdots, R_{\alpha_s})}
({\cal L}_{\alpha_1, \cdots, \alpha_s}).} 

\medskip\noindent{\bf Examples}: 

\noindent$\bullet$ For the trefoil knot $3_1$, the quantum-group invariants for
the  lowest representations are \lm:
\eqn\trefoil{
\eqalign{
W_{\tableau{1}} =& {1 \over q^{1\over 2}
 -q^{-{1\over 2}} }(-2  \lambda^{1\over2} + 3  \lambda^{3\over2} 
- \lambda^{5\over2}) + (q^{1\over 2} - q^{-{1\over 2}})(-\lambda^{1\over2} +  
\lambda^{3\over2}),\cr
W_{\tableau{2}}=&
 {(\lambda-1)(\lambda q-1)  \over
 \lambda ( q^{{1\over 2}} -
 q^{-{1\over 2}}) ^2 \,( 1 + q)} \cr
& \,\,\,\,\,\, \times 
\Bigl((\lambda q^{-1})^2
( 1 - {\lambda}q^2 +q^3  -\lambda q^3+ q^4 -\lambda q^5 
 + \lambda^2 q^5 +q^6 -\lambda q^6) \Bigr),\cr 
W_{\tableau{1 1}}=&
{(\lambda-1)(\lambda-q) \over \lambda  (  q^{{1\over 2}} -
 q^{-{1\over 2}}) ^2\,\,( 1 + q )  } \cr
& \,\,\,\,\,\, \times 
\Bigl( (\lambda q^{-2})^2 (1 -{\lambda} - \lambda q  + {{{\lambda}}^2} q
+ q^2 +q^3 - {\lambda}q^3 -{{\lambda}}\,q^4 +q^6 ) \Bigr).\cr}}

\noindent$\bullet$ For the Hopf link, one finds:
\eqn\hopfl{
W_{(\tableau{1},\tableau{1})}=\biggl( {\lambda^{1\over2}-\lambda^{-
{1\over2}} 
\over q^{1\over2} -q^{-{1\over 2}}}\biggr)^2-\lambda^{-1}(\lambda-1).}

\newsec{Reformulated polynomial invariants and a conjecture 
on their structure}

In this section we will introduce a generating functional for the quantum-group
polynomial invariants which will lead to their reformulation in terms of new
polynomials. In addition, we will state a conjecture on their structure.

\subsec{Generating functional for quantum-group invariants}

To state the conjecture about the structure of the quantum-group polynomial
invariants of links, it is useful to package these  invariants into a
generating functional. The quantum-group  invariants are labeled by
representations of $SU(N)$, which  can be also regarded as representations of
the permutation  group. It is useful to introduce a related set of
quantum-group  invariants which are labeled by conjugacy classes of the 
permutation group. We will specify these conjugacy classes by vectors 
$\vec k=(k_1, k_2, \cdots)$ with,
\eqn\sumas{
\ell=\sum_j j \, k_j,  \,\,\,\,\,\ |\vec k| =\sum_j k_j, }
and $\ell$ is then the order of the permutation group.  
The conjugacy class associated to a vector $\vec k$ will be 
denoted by $C(\vec k)$, and it has $k_j$ cycles of length $j$. 
The number of elements in the class, denoted by $|C(\vec k)|$,  
is given by the formula,
\eqn\conj{
|C(\vec k)|={\ell! \over \prod k_j! \prod j^{k_j}}.}
Given a link ${\cal L}$ with $L$ components, and given 
$L$ vectors $\vec k^{(1)}, \cdots, \vec k^{(L)}$, we define the following 
linear combination of quantum-group invariants:
\eqn\kbase{
W_{(\vec k^{(1)}, \cdots, \vec k^{(L)})}=\sum_{\{  R_{\alpha} \} } 
\prod_{\alpha=1}^L \chi_{R_\alpha}(C(\vec k^{(\alpha)})) W_{(R_1, \cdots,
R_L)},} where $\chi_{R_\alpha}$ are characters of the symmetric group
$S_{\ell_{\alpha}}$, and $\ell_{\alpha}=\sum_j j k_j^{(\alpha)}$ 
equals the number of boxes in the Young tableau of $R_{\alpha}$. If $\vec
k^{(\alpha)}$ is the zero vector for some $\alpha$, then  a 
$R_{\alpha}$ will be the trivial representation.

We now introduce generic $SU(N)$ elements,
$V_{\alpha}$, $\alpha=1, \cdots, L$, 
one for each component of the link, and denote:
\eqn\kabasev{
\Upsilon_{\vec k}(V_\alpha)=\prod_j \bigl( {\rm Tr}
\, V_\alpha^j\bigr)^{k_j}, \,\,\,\,\,\,\,\,\,\,\,\,\,\,\,\, \alpha=1,
\cdots, L.} 
If $x_j^{(\alpha)}$, $j=1, \cdots, N$ are the eigenvalues of $V_{\alpha}$, 
then \kabasev\ are symmetric polynomials in these eigenvalues (the Newton 
polynomials) and they provide a basis for ${\bf Q}[x_1^{(\alpha)}, \cdots, 
x_N^{(\alpha)}]^{S_N}$, the ring of symmetric functions in $x_1^{(\alpha)},
\cdots,  x_N^{(\alpha)}$ with rational coefficients. It is convenient in fact to 
consider $V_{\alpha} \in SU(\infty)$, for all $\alpha=1, \cdots, L$, and 
correspondingly to consider the ring of symmetric functions (in fact power 
series) in infinitely many variables, 
$x_1^{(\alpha)}, x_2^{(\alpha)}, \cdots$, which will be denoted by 
$\Lambda^{(\alpha)}$. Armed with this machinery we make 
the following definition.

\medskip\noindent{\bf Definition}: Let ${\cal L}$ be a link with $L$
components. The {\it generating functional of quantum-group
polynomial invariants} is formally defined as  follows:
\eqn\zv{
Z(V_1, \cdots, V_L)= 1 + \sum_{\{ \vec k^{(\alpha)} \}} W_{(\vec k^{(1)}, 
\cdots, \vec k^{(L)})}\prod_{\alpha=1}^L {|C(\vec k^{(\alpha)})|
\over \ell_{\alpha}!} \Upsilon_{\vec k^{(\alpha)}}(V_{\alpha}),}
where the sum is over all vectors $\vec k^{(\alpha)}$ such that 
$\sum_{\alpha=1}^L |\vec k^{(\alpha)}|>0$. Notice that, if $\vec k^{(\alpha)}$ 
is the zero vector for some $\alpha$, the invariant defined by \kbase\ will 
involve quantum group invariants with trivial 
representations. Therefore, the generating functional \zv\ 
of the quantum group invariants of a link also includes all the quantum group 
invariants of all its sublinks. Also notice that $Z(V_1, \cdots, V_L)$ 
can be regarded as as element in the ring $\Lambda^L 
(q^{\pm 1/2}, \lambda^{\pm 1/2})$ of rational functions in 
$q^{\pm 1/2}, \lambda^{\pm 1/2}$ with coefficients in $\Lambda^L =
\otimes_{\alpha=1}^L \Lambda^{(\alpha)}$.  

We will also use the logarithm of this generating functional 
(also known as {\it free energy}), which can be written as:
\eqn\fe{
F(V_1, \cdots, V_L)= \log \,Z(V_1, \cdots, V_L)=
 \sum_{\{ \vec k^{\alpha} \} } W^{(c)}_{(\vec k^{(1)}, 
\cdots, \vec k^{(L)})}\prod_{\alpha=1}^L {|C(\vec k^{(\alpha)})|
\over \ell_{\alpha}!} \Upsilon_{\vec k^{(\alpha)}}(V_{\alpha}),} 
and defines the ``connected'' invariants $W^{(c)}_{(\vec k^{(1)}, 
\cdots, \vec k^{(L)})}$. 
 
\medskip\noindent{\bf Examples}: In the case of a knot, one has for the
simplest cases,
\eqn\exem{
\eqalign{
W^{(c)}_{(1,0, \dots)}=&W_{(1,0,\cdots)}=W_{\tableau{1}},\cr
W^{(c)}_{(2,0, \dots)}=&W_{(2,0,\cdots)}-W^2_{(1,0,\cdots)} = 
W_{\tableau{2}} + W_{\tableau{1 1}}-W^2_{\tableau{1}},\cr}}
while for a 2-component link ${\cal L}$ with components ${\cal K}_1$ and
${\cal K}_2$,
\eqn\exemlink{
\eqalign{
W^{(c)}_{((1,0, \dots),(1,0,
\dots))}({\cal L})=&W_{((1,0, \dots),(1,0,
\dots))}({\cal L})-W_{(1,0,
\dots)}({\cal K}_1)W_{(1,0,\dots)}({\cal K}_2)\cr=&
W_{(\tableau{1},\tableau{1})}({\cal
L})-W_{\tableau{1}}({\cal K}_1) W_{\tableau{1}}({\cal K}_2).\cr}}

\subsec{Reformulated polynomial invariants}

Our next goal is to reformulate the quantum-group polynomial invariants. Let
us begin introducing rational functions
$f_{(R_1, \cdots, R_L)} (q,\lambda)$, labeled by  representations of $SU(N)$,
as follows:
\eqn\ovc{
F(V_1, \cdots, V_L)=\sum_{d=1}^{\infty} \sum_{ \{ R_{\alpha}\}  } 
{1 \over d} f_{(R_1, \cdots, R_L)}(q^d, \lambda^d) \prod_{\alpha=1}^L 
{\rm Tr}_{R_{\alpha}}V_{\alpha}^d.}
Using Frobenius formula and \fe, the relation \ovc\ is equivalent to
the  following equation:
\eqn\master{
W_{(\vec k^{(1)}, \cdots,\vec k^{(L)}) }^{(c)}=\sum_{d|\vec k^{(\alpha)}} 
d^{\sum_{\alpha} |\vec k^{(\alpha)}|-1} \sum_{ \{ R_{\alpha}\} }
\prod_{\alpha=1}^L \chi_{R_{\alpha}}(C(\vec k_{1/d}^{(\alpha)})) 
f_{(R_1, \cdots, R_L)}(q^d, \lambda^d).} 
In this equation, the vector $\vec k_{1/d}$ is defined as 
follows. Fix a vector $\vec k$, and 
consider all the positive integers $d$ that satisfy the following condition:
$d|j$ for every $j$ with $k_j \not=0$. When this happens, 
we will say that ``$d$ divides 
$\vec k$", and we will denote this as $d|\vec k$.  
We can then define the vector $\vec k_{1/d}$ whose components are:
\eqn\shifted{
(\vec k_{1/d})_i=k_{di}.}
The vectors which satisfy the above condition and are
 ``divisible by $d$'' have the structure 
$(0,\dots, 0,k_d, 0,\dots,0, k_{2d},0,\dots)$, and the vector 
$\vec k_{1/d}$ is then
given by $(k_d,k_{2d},\dots)$. In \master, the integer $d$ has to 
divide all the vectors $\vec k^{(\alpha)}$, $\alpha=1, \cdots, L$. 

It is easy to show that 
the rational functions $f_{(R_1, \cdots, R_L)}$ 
are uniquely determined in terms of the connected invariants $W_{(\vec
k^{(1)}, \cdots,\vec k^{(L)}) }^{(c)}$ by equation \master. Moreover, one can 
invert this equation and write an explicit formula for 
$f_{(R_1, \cdots, R_L)}$ in terms of the quantum group invariants \foot{We 
are grateful to Ezra Getzler for explaining to us how to find this 
explicit expression. A similar inversion formula has been 
obtained for the Gopakumar-Vafa invariants in 
\bp.}. Following \gk, if $F\in \Lambda^L (q^{\pm 1/2}, \lambda^{\pm 1/2})$, 
we can define an operation $\psi_d$ on this ring by 
\eqn\opera{
\psi_d \circ F(q,\lambda; \Upsilon_{\vec k^{(\alpha)}}(V_{\alpha}))= 
F(q^d, \lambda^d;\Upsilon_{\vec k^{(\alpha)}}(V^d_{\alpha})).}
We can also define the so-called plethystic exponential \gk, 
\eqn\pletexp{
{\rm Exp}(F)= \exp\Bigl( \sum_{d=1}^{\infty}{\psi_d \over d}\Bigr)\circ F.}
This exponential has an inverse given by \gk:
\eqn\logplet{
{\rm Log} (F)= \sum_{d=1}^{\infty} {\mu (d) \over d} \log (\psi_d \circ F),}
where $\mu(d)$ is the Moebius function. 
Recall that this function is zero if $d$ is not square-free, and  
it is $(-1)^s$ otherwise, where $s$ is the number of primes in the 
decomposition of $d$. 

Using these results, we can write an 
explicit formula for the $f_{(R_1, \cdots, R_L)}(q, \lambda)$. 
From \opera\ and \pletexp\ it follows that
\eqn\zplet{
Z(V_1, \cdots, V_L)= {\rm Exp} \biggl( \sum_{ \{ R_{\alpha} \}} 
f_{(R_1, \cdots, R_L)}(q, \lambda) \prod_{\alpha=1}^L {\rm Tr}_{R_{\alpha}}V_{\alpha}
\biggr).}
Using now the inverse of the plethystic exponential and the Frobenius 
formula, we finally obtain:
\eqn\explinverlink{
\eqalign{
f_{(R_1, \cdots, R_L)}(q, \lambda) =&
 \sum_{d, m=1}^{\infty} (-1)^{m-1} {\mu 
(d) \over d m} 
\sum_{ \{ \vec k^{(\alpha \, j)}, R_{\alpha\, j} \} }
\prod_{\alpha=1}^L 
 \chi_{R_{\alpha}} \biggl( 
C\biggl( (\sum_{j=1}^m \vec k^{(\alpha \, j)})_d\biggr)\biggr)\cr
& \times \prod_{j=1}^m {|C(\vec k^{(\alpha \, j)})| \over \ell_{\alpha \, j}!}
  \chi_{R_{\alpha\, j}}(C(\vec k^{(\alpha \, j)})) 
W_{(R_{1j}, \cdots, R_{Lj})}(q^d, \lambda^d).}}
The second sum runs over all vectors ${\vec k}^{(\alpha \, j)}$, 
with $\alpha=1, \cdots, L$ and $j=1, \cdots, m$,  
such that $\sum_{\alpha =1}^L |\vec k^{(\alpha\, j)}| >0$ 
for any $j$, and over representations $R_{\alpha \, j}$.  
In \explinverlink, the vector $\vec k_d$ is defined 
as follows: $(\vec k_d)_{di}=k_i$ and has zero 
entries for the other components. Therefore, if $\vec k= (k_1, k_2, \cdots)$, 
then
\eqn\multd{
\vec k_d =(0, \cdots, 0,k_1,0,\cdots, 0, k_2, 0,\cdots),}
where $k_1$ is in the $d$-th entry, $k_2$ in the $2d$-th entry, and so 
on. Notice that $(\vec k_d)_{1/d}=\vec k$, so this is the inverse 
operation to \shifted. Using \explinverlink, 
it is easy to show that the $f_{(R_1, \cdots, R_L)}$ are equal 
to the quantum-group invariants $W_{(R_1, \cdots,R_L)}$, plus some 
extra terms that involve the $W_{(R'_1, \cdots,R'_L)}$ 
with a lower total number of boxes $\sum_{\alpha} \ell'_{\alpha}$. These
functions $f_{(R_1, \cdots, R_L)}$ are indeed the reformulated 
polynomial invariants.

\medskip\noindent{\bf Definition}: The reformulated quantum-group polynomial
invariants are the functions $f_{(R_1, \cdots, R_L)}$ entering \ovc. 
They can be expressed in terms of quantum group invariants through 
\explinverlink.

\medskip\noindent{\bf Examples}: In the case of knots, one has, for 
representations of up to three boxes:
\eqn\exov{
\eqalign{
f_{\tableau{1}}(q,\lambda)=&W_{\tableau{1}}(q, \lambda), 
\cr
f_{\tableau{2}}(q,\lambda)=&W_{\tableau{2}}(q,\lambda)
-{1\over 2}\bigl( W_{\tableau{1}}(q,\lambda)^2+ 
W_{\tableau{1}}(q^2,\lambda^2)
\bigr),\cr
f_{\tableau{1 1}}(q, \lambda)=&W_{\tableau{1 1}}(q,\lambda)
-{1\over 2}\bigl( W_{\tableau{1}}(q,\lambda)^2-
 W_{\tableau{1}}(q^2,\lambda^2) 
\bigr),\cr 
f_{\tableau{3}}(q, \lambda) =&
W_{\tableau{3}}(q, \lambda)-W_{\tableau{1}}(q, \lambda)
W_{\tableau{2}}(q, \lambda)+{1 \over 3} (W_{\tableau{1}}(q, \lambda)^3 
-W_{\tableau{1}}(q^3, \lambda^3)),\cr
f_{\tableau{2 1}}(q, \lambda)=& W_{\tableau{2 1}}(q, \lambda) - 
W_{\tableau{1}}(q, \lambda) 
(W_{\tableau{2}}(q, \lambda) + W_{\tableau{1 1}}(q, \lambda))+ 
{2 \over 3} W_{\tableau{1}}(q, \lambda)^3+
{1 \over 3} W_{\tableau{1}}(q^3, \lambda^3),\cr
f_{\tableau{1 1 1}}(q, \lambda)=& W_{\tableau{1 1 1}}(q, \lambda) - 
W_{\tableau{1}}(q, \lambda) 
W_{\tableau{1 1}}(q, \lambda) + 
{1 \over 3} (W_{\tableau{1}}(q, \lambda)^3-
 W_{\tableau{1}}(q^3, \lambda^3)).
\cr} }
On the other hand, for a link ${\cal L}$ with two components, ${\cal K}_1$ and
${\cal K}_2$,  one finds:
\eqn\linkcompo{
\eqalign{
f_{(\tableau{1}, \tableau{1})}({\cal L})= &
 W_{(\tableau{1}, \tableau{1})}({\cal L}) - 
W_{\tableau{1}}({\cal K}_1) W_{\tableau{1}}({\cal K}_2),\cr
f_{(\tableau{2}, \tableau{1})}({\cal L})=&
W_{(\tableau{2}, \tableau{1})}({\cal L})-
W_{(\tableau{1}, \tableau{1})}({\cal L})W_{\tableau{1}}({\cal K}_1) 
-W_{\tableau{2}}({\cal K}_1)W_{\tableau{1}}({\cal K}_2)+ 
W_{\tableau{1}}({\cal K}_1)^2 W_{\tableau{1}}({\cal K}_2),\cr 
f_{(\tableau{1 1}, \tableau{1})}({\cal L})=&
W_{(\tableau{1 1}, \tableau{1})}({\cal L})-
W_{(\tableau{1}, \tableau{1})}({\cal L})W_{\tableau{1}}({\cal K}_1) 
-W_{\tableau{1 1}}({\cal K}_1)W_{\tableau{1}}({\cal K}_2)+ 
W_{\tableau{1}}({\cal K}_1)^2 W_{\tableau{1}}({\cal K}_2).\cr} 
}
\medskip
\noindent
{\bf Remark}: The functions $f_{(R_1, \cdots, R_L)}$ were introduced 
by Ooguri and Vafa \ov\ through the relation \ovc. A recursive procedure 
to obtain these functions in terms of quantum-group polynomial invariants 
was spelled out in detail in \lm\lmv.   

\subsec{A conjecture on the structure of the reformulated polynomial 
invariants}

In this subsection we present a conjecture on the algebraic structure
of $f_{(R_1, \cdots, R_L)}$, which 
in turn implies a structure result for the $W_{(R_1, 
\cdots, R_L)}$. To state it we need to introduce the 
Clebsch-Gordan coefficients of the symmetric group. 
These are given by,
\eqn\cgcoef{
C_{R R' R''}=\sum_{\vec k} {|C(\vec k)| \over \ell!} \chi_R (C(\vec k)) 
\chi_{R'} (C(\vec k)) \chi_{R''} (C(\vec k)).}
Finally, we also need to introduce the monomial $S_R(q)$, defined as
follows.  If $R$ is a hook representation, {\it i.e.}, a representation whose 
Young tableau is of the form,
\eqn\hook{
\tableau{10 1 1 1 1 1}
}
with $\ell$ boxes in total and $\ell-d$ boxes in the first row, then
\eqn\mono{
S_R(q)=(-1)^d q^{-{\ell-1\over 2} +d}, }
and $S_R(q)=0$ otherwise. Now we are ready to formulate the conjecture.

\medskip\noindent{\bf Conjecture }: Given a link ${\cal L}$, the
reformulated quantum-group polynomial invariants,
$f_{(R_1, \cdots, R_L)}(q, \lambda)$, have the following 
structure: 
\eqn\frlinks{
\eqalign{
&f_{(R_1, \cdots, R_L)}(q, \lambda)=\cr & 
(q^{{1\over 2}}-q^{-{1 \over 2}})^{L-2}
\sum_{g\ge 0} \sum_{Q}
\sum_{\{R'_{\alpha}, R_{\alpha}''\} } \Bigl( \prod_{\alpha=1}^L
C_{R_{\alpha}\,R_{\alpha}'\,R_{\alpha}''}S_{R_{\alpha}'}(q)\Bigr) 
N_{(R_1'', \cdots, R_L''),g,Q} 
 (q^{{1\over 2}}-q^{-{1\over 2}})^{2g}\lambda^Q,\cr}}
where $N_{(R_1, \cdots, R_L),g,Q}$ are {\it integer} numbers\foot{These 
integers differ in a
sign  from the integers introduced in \lmv. More precisely, the 
$N_{(R_1, \cdots, R_L),g,Q}$ are $(-1)^{L-2}$ times the integers 
denoted by ${\widehat N}_{(R_1, \cdots, R_L),g,Q}$ in \lmv.}, and $Q$ are
either all integers or all semi-integers.

\medskip\noindent{\bf Remark}: the conjecture \frlinks\ was proposed in 
\lmv, and refines and generalizes a previous conjecture by 
Ooguri and Vafa. It can be 
regarded as a definition of the integer invariants 
$N_{(R_1, \cdots, R_L),g,Q}$. The fact that they 
can be extracted from the $W_{(R_1, \cdots, R_L)}(q, \lambda)$ 
in the way described is far from obvious. For example, according 
to the conjecture, $f_{(R_1, \cdots, R_L)}$ must be a 
polynomial in $q^{\pm {1 \over 2}}, \lambda^{\pm {1 \over 2}}$ 
with integer coefficients, times 
$(q^{1 \over 2} -q ^{-{1 \over 2}})^{L-2}$.

\medskip\noindent{\bf Examples}: 

\noindent$\bullet$ For a knot ${\cal K}$ with HOMFLY
polynomial $P_{\cal K}$, one has,
\eqn\homk{
P_{\cal K} (q, \lambda)=\sum_{g, j} a_{g,j}
(q^{1 \over 2} -q ^{-{1 \over 2}})^{2g} \lambda^j,
}
and the integer invariants are:
\eqn\inth{
N_{\tableau{1}, g, j+1/2}=a_{g,j+1}-a_{g, j}.
}

\noindent$\bullet$ For the trefoil knot, after using the known quantum-group 
invariants, one finds \lm:
\eqn\ftrefoil{
\eqalign{
f_{\tableau{2}}(q,\lambda)
&={ q^{-{1\over 2}}{\lambda}( {\lambda}-1) ^2 \,\,( 1 + {q^2}) \,
     ( q + {{{\lambda}}^2}\,q - {\lambda}\,( 1 + {q^2} ))
\over q^{{1\over 2}} - q^{-{1\over 2}} },  
\cr
f_{\tableau{1 1}}(q,\lambda)& = - {1 \over q^3}f_{\tableau{2}}(q,\lambda).\cr}
}
The corresponding integer invariants are listed in the following tables:
 
\medskip

{\vbox{\ninepoint{
$$
\vbox{\offinterlineskip\tabskip=0pt
\halign{\strut
\vrule#
&
&\hfil ~$#$
&\hfil ~$#$
&\hfil ~$#$
&\hfil ~$#$
&\hfil ~$#$
&\hfil ~$#$
&\vrule#\cr
\noalign{\hrule}
&g
&Q=1
&2
&3
&4
&5
&
\cr
\noalign{\hrule}
&0
&-2
&8
&-12
&8
&-2
&\cr
&1
&-1
&6
&-10
&6
&-1
&\cr
&2
&0
&1
&-2
&1
&0
&\cr
\noalign{\hrule}}\hrule}$$}
\vskip - 7 mm
\centerline{{\bf Table 1:} The integers $N_{\tableau{2},g,Q}$
for 
the trefoil knot.}
\vskip7pt}
\noindent
\smallskip

{\vbox{\ninepoint{
$$
\vbox{\offinterlineskip\tabskip=0pt
\halign{\strut
\vrule#
&
&\hfil ~$#$
&\hfil ~$#$
&\hfil ~$#$
&\hfil ~$#$
&\hfil ~$#$
&\hfil ~$#$
&\vrule#\cr
\noalign{\hrule}
&g
&Q=1
&2
&3
&4
&5
&
\cr
\noalign{\hrule}
&0
&-4
&16
&-24
&16
&-4
&\cr
&1
&-4
&20
&-32
&20
&-4
&\cr
&2
&-1
&8
&-14
&8
&-1
&\cr
&3
&0
&1
&-2
&1
&0
&\cr
\noalign{\hrule}}\hrule}$$}
\vskip - 7 mm
\centerline{{\bf Table 2:} The integers $N_{\tableau{1 1},g,Q}$ 
for 
the trefoil knot.}
\vskip7pt}
\noindent

\subsec{Some consequences for the algebraic structure 
of the HOMFLY polynomial of links}

In this subsection, we will explore some of the 
consequences of \frlinks\ for the algebraic structure 
of the HOMFLY polynomial of links.
To do this, it is convenient to introduce some notation.  
Let us consider a link ${\cal L}$ of $L$ components ${\cal K}_{\alpha}$, 
$\alpha=1, \cdots, L$. For simplicity we will denote 
\eqn\sim{
W_{(\tableau{1}, 
\cdots, \tableau{1})}({\cal L})=W_{((1,0,\cdots), \cdots, (1,0,\cdots))} =
W({\cal L}).} To write the ``connected'' invariant defined by \fe, we 
have to take into account the invariants of all the sublinks of ${\cal L}$. 
In section 2, given a subset $\{\alpha_1, 
\cdots, \alpha_s\} \subset \{1, \cdots, L\}$, we defined  
${\cal L}_{\alpha_1, \cdots, \alpha_s}$ as the sublink of $s$ components 
which is obtained from the link ${\cal
L}$ by keeping the components ${\cal K}_{\alpha_i}$, $i=1, \cdots, s$, and 
by ``deleting'' the remaining $L-s$ components. One can easily see from the 
definition of the ``connected'' invariants that  
$W^{(c)}({\cal L})$ is given by 
the original invariant \sim\ plus some extra terms involving 
products of invariants for sublinks: 
\eqn\conn{
W^{(c)}({\cal L})= W({\cal L})
-\sum_{{\alpha_L}=1}^L
 W({\cal L}_{\alpha_1, \cdots, \alpha_{L-1}})\, 
W({\cal L}_{\alpha_L}) 
+ \cdots.}
For example, for a link of two components ${\cal K}_1$, and ${\cal K}_2$ 
one simply has:
\eqn\twoc{
W^{(c)}({\cal L})=W({\cal L})-W({\cal K}_1) W({\cal K}_2).}

Using \master, we see that the conjecture \frlinks\ states that, 
\eqn\linkst{
W^{(c)}({\cal L})=(q^{{1 \over 2}} -q^{-{1\over 2}})^{L-2} 
\sum_Q \sum_{g\ge 0} N_{(\tableau{1}, \cdots, \tableau{1}),g,Q} 
\lambda^Q (q^{{1\over 2}}-q^{-{1\over 2}})^{2g}.}
To analyze \linkst, we will first consider the simple case of a 
link of two components.
Using \twoc\ and \rel, we find that the HOMFLY polynomial of the link has the 
following structure: 
\eqn\twocstr{
P_{\cal L}(q, \lambda)=\sum_{g\ge 0} p^{\cal L}_{2g-1} 
(\lambda)(q^{1\over 2}-q^{-{1\over 2}})^{2g-1},} 
{\it i.e.} the lowest power of $q^{{1\over 2}}-q^{-{1\over 2}}$ is $-1$, 
and the powers are congruent to $-1$ mod 2. Moreover, if we denote the 
HOMFLY polynomial of the component knots by,
\eqn\compknot{
P_{{\cal K}_{\alpha}}(\lambda,q)=\sum_{g\ge 0}p^{{\cal K}_\alpha}_{2g}
(\lambda)
(q^{{1\over 2}}-q^{-{1\over 2}})^{2g},}
for $\alpha=1,2$, we find,
\eqn\lmthtwo{
p_{-1}^{\cal L}(\lambda)=\lambda^{-{\rm lk}({\cal L})}(\lambda^{1\over2}-
\lambda^{-{1\over2}})p_0^{{\cal K}_1}(\lambda)p_0^{{\cal K}_2}(\lambda).} 
The last equation comes from the requirement that there are no powers of 
$(q^{{1\over 2}}-q^{-{1\over 2}})^{-2}$ in  
$W^{(c)}({\cal L})$. The results \twocstr\ and \lmthtwo\ 
capture completely the algebraic structure of the HOMFLY polynomial of a 
two-component link, and reproduce the results of Lickorish 
and Millett \lickm.

We can generalize the above results for links with an arbitrary number 
of components $L$. By induction on the number of components, and using 
\conn\ and \linkst, it is easy to prove that the HOMFLY polynomial of the
link has  the following structure:
\eqn\link{
P_{\cal L}(q,\lambda)=
\sum_{g\ge0} p^{\cal L}_{2g+1-L} 
(\lambda) (q^{{1\over 2}}-q^{-{1\over 2}})^{2g+1-L},}
{\it i.e.} the lowest power of $q^{{1\over 2}}-q^{-{1\over 2}}$ is $1-L$. 
This has been proved in \lickm. 
Due to the relation \rel, it is convenient to 
introduce the following polynomials in $\lambda$:
\eqn\newpol{
{\widetilde p}^{{\cal L}_{\alpha_1, \cdots, \alpha_s}}_k 
(\lambda)={\lambda}^{{\rm lk}({\cal L}_{\alpha_1, \cdots, \alpha_s})} 
p^{{\cal L}_{\alpha_1, \cdots, \alpha_s}}_k 
(\lambda).} 
Finally, we will write,
\eqn\connagain{
W^{(c)}({\cal L})=\biggl( {\lambda^{1\over2}-
\lambda^{-{1\over2}}\over q^{1\over2}-q^{-{1\over2}}}\biggr)\sum_{g\ge0} 
{\widetilde p}^{(c),{\cal L}}_{2g+1-L} 
(\lambda)(q^{1\over2}-q^{-{1\over2}})^{2g+1-L}.}  
The conjecture \linkst\ then states that,
\eqn\strhom{
{\widetilde p}_{1-L}^{(c),{\cal L}}(\lambda)=
{\widetilde p}_{3-L}^{(c),{\cal L}}(\lambda)=
\cdots ={\widetilde p}_{L-3}^{(c),{\cal L}}(\lambda)=0.} This implies, 
in particular, that the polynomials $p^{\cal L}_k(\lambda)$ 
of the HOMFLY polynomial of a link, 
for $k=1-L, 3-L, \cdots, L-3$, are completely determined by the HOMFLY 
polynomial of its sublinks.  As a first consequence of \strhom, it is easy to
show the  following proposition.

\medskip\noindent{\bf Proposition} (Lickorish and Millett \lickm). The
polynomial in
$\lambda$, 
$p_{1-L}(\lambda)$, in the HOMFLY polynomial of a link \link\ is 
given by
\eqn\lickmth{
p_{1-L}^{\cal L}(\lambda)=\lambda^{-{\rm lk}({\cal L})}(\lambda^{1\over2}- 
\lambda^{-{1\over2}})^{L-1} \prod_{\alpha=1}^L p_0^{{\cal
K}_{\alpha}}(\lambda).}

This result is a consequence of 
${\widetilde p}_{1-L}^{(c),{\cal L}}(\lambda)=0$, and it is 
easily proven by induction on the number of components of the link: 
since ${\widetilde p}_{1-L}^{(c), {\cal L}}(\lambda)=0$, one can extract
the coefficient of the lowest power of $q^{1\over2}-q^{-{1\over2}}$ in 
$W({\cal L})$ from the terms in the expansion of $W^{(c)}({\cal L})$ that 
only involve products of invariants of 
sublinks. One sees immediately that 
the relevant part of these invariants is again the coefficient of 
the lowest power of $q^{1\over2}-q^{-{1\over2}}$. But because of the 
induction hypothesis, these in turn can be evaluated by factorization 
into their knots. This means that the coefficient 
${\widetilde p}_{1-L}^{\cal L}$  can be evaluated from $\prod_{\alpha=1}^L 
W({\cal K}_{\alpha})$, and this proves \lickmth. 

Notice that \lickmth\ is just the simplest 
consequence of \strhom, which gives 
much more relations. For example, for links with $L=3$, 
the equality ${\widetilde p}_{0}^{(c), {\cal L}}(\lambda)=0$ implies that
\eqn\lasteq{
\eqalign{
 {\widetilde p}^{\cal L}_{0}(\lambda)&=(\lambda^{1\over2}-
\lambda^{-{1\over2}})(p_0^{{\cal K}_1}(\lambda) \,
{\widetilde p}_1^{{\cal L}_{23}}(\lambda) + {\rm perms}) \cr 
&-2(\lambda^{1\over2}-\lambda^{-{1\over 2}})^2(p_2^{{\cal K}_1}(\lambda)
\, p_0^{{\cal K}_2}(\lambda)\,p_0^{{\cal K}_3}(\lambda)+ {\rm
perms}).\cr}}
For links with more components, one obtains more complicated 
equations which can be summarized as in \strhom, providing
a new set of results on the algebraic 
structure of the HOMFLY polynomial of links. 

\newsec{Topological content  of the new integer invariants}

In the previous sections, we have presented the conjecture 
on the structure of the reformulated quantum-group invariants of
knots and links, and  we have introduced a new set of integer invariants.
In this  section, we will describe the topological content of the latter. 
The starting point is the connection between quantum-group polynomial
invariants and Chern-Simons gauge theory. Quantum-group invariants can be
expressed as vacuum expectation values of Wilson loops. As in any gauge
theory these vacuum expectation values admit a large-$N$ expansion \thooft, 
which in turn can be interpreted as a string theory expansion\foot{The
large-$N$  expansion of Wilson loops in Chern-Simons gauge theory was
studied from a  field theory point of view in \cg.}. The string 
theory description of Wilson loops in Chern-Simons gauge theory 
is given in terms of a  
topological open string theory, of the kind that in the closed case leads to
Gromov-Witten invariants. The first step to provide a geometrical meaning
to the new integer invariants is to express the reformulated quantum-group
polynomial invariants in terms of Gromov-Witten invariants generalized to
the open string case.

\subsec{$1/N$ expansion and Gromov-Witten 
invariants for open strings}

The geometric interpretation first appears 
in the context of the so-called $1/N$ expansion 
of the invariants. The ``connected'' invariants that 
we introduced in \fe\ are rational 
functions of $q^{\pm {1\over 2}}$ and 
$\lambda^{\pm {1\over 2}}$. If we put $q={\rm e}^{ix}$ but 
keep $\lambda$ fixed, and we formally expand in $x$, we 
find a series with the structure:
\eqn\largen{
\Bigl( \prod _{\alpha=1}^L {|C(\vec k^{(\alpha)})| \over 
\ell_{\alpha}!}\Bigr) W^{(c)}_{(\vec k^{(1)}, \cdots, \vec k^{(L)})}
=i^{\sum_{\alpha=1}^L |\vec k^{(\alpha)}|}\sum_{g=0}^{\infty} 
x^{2g-2+\sum_{\alpha=1}^L |\vec k^{(\alpha)}|} F_{g, 
(\vec k^{(1)}, \cdots, \vec k^{(L)}) }(\lambda).}
Notice that the generating functional \fe\ can be written in terms 
of the functions $F_{g, 
(\vec k^{(1)}, \cdots, \vec k^{(L)}) }(\lambda)$ 
as,
\eqn\fef{F(V_1, \cdots, V_L)= 
 \sum_{\{ \vec k^{(\alpha)} \} }i^{\sum_{\alpha=1}^L 
|\vec k^{(\alpha)}|}\sum_{g=0}^{\infty}  
x^{2g-2+\sum_{\alpha=1}^L |\vec k^{(\alpha)}|}
F_{g,(\vec k^{(1)}, \cdots, \vec k^{(L)}) }(\lambda) 
\prod _{\alpha=1}^L\Upsilon_{\vec k^{(\alpha)}}(V_{\alpha}).}
 
\medskip\noindent{\bf Remarks}:

\noindent$\bullet$ Since we are not expanding $\lambda={\rm e}^{iNx}$, 
we are keeping the variable $t=iNx$ fixed, and therefore we 
can equivalently write the above series as a ``$1/N$ expansion'' 
by putting $x=-it/N$. The parameter $t$ is also 
called the 't Hooft parameter.

\noindent$\bullet$ The structure of the above expansion can be 
proved in the context of Chern-Simons theory by using 
standard $1/N$ analysis. We are not aware of a proof 
relying on the quantum-group definition of the invariants.

\medskip

The geometric picture for the reformulated quantum-group invariants is 
based on the proposals made in \gv\ov. Before stating it we need
to introduce some machinery.

It was conjectured in \ov\ that to every link ${\cal L}$ in ${\bf
S}^3$ one can associate a Lagrangian submanifold ${\cal C}_{\cal L}$ in
the non-compact  Calabi-Yau $X$,
\eqn\rescon{
{\cal O}(-1)\oplus{\cal O}(-1) \rightarrow {\bf P}^1, }
also called the {\it resolved conifold}. 
The assignment implies that $b_1({\cal C}_{\cal L})=L$, the number 
of components of ${\cal L}$.

 The quantities $F_{g,(\vec k^{(1)}, \cdots, \vec k^{(L)})
}(\lambda)$ are then interpreted in terms of an appropriate generalization 
of the Gromov-Witten invariants for Riemann surfaces 
with boundaries. Let $\gamma_{\alpha}$, $\alpha=1, 
\cdots, L$, be one-cycles representing a basis for 
$H_1 ({\cal C}_{\cal L}, {\bf Z})$, and let 
${\cal Q} \in H_2(X, {\cal C}_{\cal L}, {\bf Z})$ 
be a relative two-homology class ({\it i.e.}, 
a two-cycle of $X$ that ends on ${\cal C}_{\cal L}$). 
Then, one considers the
holomorphic maps $f: \Sigma_{g,h} \rightarrow X$ 
from a Riemann surface $\Sigma_{g,h}$ of genus 
$g$ and with $h =\sum_{\alpha=1}^L |\vec k^{(\alpha)}|$ 
holes, which satisfy the following conditions: 
first, $f_*[\Sigma_{g,h}]={\cal Q}$; second,  
$k^{(\alpha)}_j$ of the $h$ (oriented) boundaries of $\Sigma_{g,h}$  
map to the cycle $\gamma_{\alpha}$ with winding number $j$, 
{\it i.e.}, $f_*[C]=j [\gamma_{\alpha}]$ for $k_j^{(\alpha)}$ 
oriented circles $C$ 
of the boundary. The ``number'' of such maps 
will be denoted by $N^{\cal Q}_{g,(\vec k^{(1)}, \cdots, \vec k^{(L)})}$. 
These numbers are the open-string analog of Gromov-Witten invariants, and   
a precise definition would involve the construction of a 
compact moduli space for the maps $f$. The invariant 
$N^{\cal Q}_{g,(\vec k^{(1)}, \cdots, \vec k^{(L)})}$ would be then 
given by the degree of the virtual fundamental class of the moduli 
space, as in Gromov-Witten theory (see \kl\ls\ for more details). 

We can now describe the geometrical content of the coefficients in the 
$1/N$ expansion. Given a link ${\cal L}$,  the 
functions $F_{g,(\vec k^{(1)}, \cdots, \vec k^{(L)})
}(\lambda)$ appearing in the expansion \largen\ are expressed in terms of the 
Gromov-Witten invariants for open strings in the following way:
\eqn\fgk{
F_{g, (\vec k^{(1)}, \cdots, \vec k^{(L)}) }(\lambda)=\sum_{{\cal Q}} 
N^{\cal Q}_{g, (\vec k^{(1)}, \cdots, \vec k^{(L)})}{\rm e}^{\int_{\cal Q} 
\omega},}
where $\omega$ is the K\"ahler class of the Calabi-Yau manifold $X$ and 
$\lambda={\rm e}^t$, with 
\eqn\kalpar{
t=\int_{ {\bf P}^1 }\omega.} 
For any ${\cal Q}$, one can always write $\int_{\cal Q}\omega =Q t$, where 
$Q$ is in general a half-integer number, and therefore 
$F_{g, (\vec k^{(1)}, \cdots, \vec k^{(L)}) }(\lambda)$ is a polynomial 
in $\lambda^{\pm {1\over 2}}$ with rational coefficients. 

\medskip\noindent{\bf Example}: In the case of the trivial knot or unknot, it
is  easy to compute the functions $F_{g, \vec k }(\lambda)$ from the 
quantum-group invariants. These functions are nonvanishing only for 
vectors of the form $\vec k=(0,\cdots, 0,1, 0, \cdots, 0)$ with 
the nonzero entry in the $d$-th position, and:
\eqn\ogw{
 F_{g, (0,\cdots, 0,1, 0, \cdots, 0)} (\lambda) ={ (1-2^{1-2g})|B_{2g}|
\over  (2g)! }d^{2g-2} ({\lambda}^{d\over2} - {\lambda}^{-{d\over 2}}).}
According to \ov, the Lagrangian submanifold ${\cal C}_{\cal K}$ 
associated to the unknot is a sphere bundle over the equator of ${\bf P}^1$, 
and the non-trivial one-cycle of ${\cal C}_{\cal K}$ is precisely this 
equator. The relative homology $H_2(X, {\cal C}_{\cal K}, {\bf Z})$ 
has two primitive generators $[N]$, $[S]$, 
corresponding to the northern and southern hemisphere 
of ${\bf P}^1$. The terms with ${\lambda}^{\pm {d \over 2}}$ correspond 
to holomorphic maps satisfying $f_* [\Sigma_{g,1}]=d[N]$ or 
$f_* [\Sigma_{g,1}]=d[S]$, respectively. The expression \ogw\
was first obtained in \ov\ by using the above conjectured 
relation with Chern-Simons theory, and it has been 
computed directly in Gromov-Witten theory in \kl\ls.
\medskip\noindent{\bf Remark}: As we said above, 
the correspondence that associates 
a Lagrangian submanifold ${\cal C}_{\cal L}$ to a ${\cal L}$ 
is conjectural, and so far there is no well-defined procedure to construct 
${\cal C}_{\cal L}$ once ${\cal L}$ is given. In \ov, Ooguri and Vafa 
showed that given a link ${\cal L}$ in ${\bf S}^3$, one can canonically 
associate to it a Lagrangian submanifold 
${\widehat {\cal C}}_{\cal L}$ in $T^*{\bf 
S}^3$. ${\cal C}_{\cal L}$ should be obtained from 
${\widehat {\cal C}}_{\cal L}$ after a ``conifold transition'' from the 
$T^*{\bf S}^3$ geometry (the deformed conifold) to the resolved conifold. 
There is however a proposal for ${\cal C}_{\cal L}$ in \lmv\ for a class 
of torus knots and links.

\subsec{The integer invariants}

To describe the geometrical content of the new integer invariants 
introduced in \frlinks, $N_{(R_1, \cdots, R_L),g,Q}$,   
we have to resum the Gromov-Witten invariants for 
open strings. It was shown in \ov\lmv\ that 
the generating functional $F(V_1, \cdots, V_L)$ can be 
written as:
\eqn\multi{
\eqalign{
 F(V_1, \cdots, V_L)= & \sum_{g, Q} \sum_{d >0} 
\sum_{\{ R_{\alpha}, R_{\alpha}', R_{\alpha}'' 
\}} 
N_{(R_1, \cdots, R_L),g,Q}\,  {1 \over d} 
 \bigl( 2i \sin (dx/2)\bigr)^{2g+L-2} \cr 
& \,\,\,\,\,\,\ \times \Bigl( \prod_{\alpha=1}^L
 C_{R_{\alpha} R_{\alpha}'R_{\alpha}''}
S_{R_{\alpha}''}({\rm e}^{d i x})\, {\rm Tr}_{R_{\alpha}} V_{\alpha}^d \Bigr)
\lambda^{d Q}.\cr}} 
The conjecture \frlinks\ follows from this equation. 
From \multi\ one can also extract the expression of the 
$F_{g, (\vec k^{(1)}, \cdots, \vec k^{(L)}) }(\lambda)$ in terms of 
integer invariants $N_{(R_1, \cdots, R_L),g,Q}$, by simply combining 
\master, \frlinks, and \largen.
 
In \lmv\ a geometric interpretation of these integer invariants 
was given in terms of the Calabi-Yau geometry 
described in the previous subsection. Let $R_1,\cdots, R_L$ be 
representations of $SU(N)$, 
where $R_{\alpha}$ has $\ell_{\alpha}$ boxes, and let ${\cal C}_{\cal L}$ 
be the Lagrangian submanifold associated to the link ${\cal L}$ with 
$L$ components. Let us denote by ${\cal M}_{g,\ell,Q}$ the 
moduli space of Riemann surfaces of genus $g$ and $\ell$ holes 
embedded in the resolved conifold, where $\ell =\sum_{\alpha=1}^L 
\ell_{\alpha}$. 
The embedding is such that $\ell_{\alpha}$ holes end on the the non-trivial 
cycles $\gamma_{\alpha}$, for $\alpha=1, \cdots, L$, and the relative 
class $H_2(X, {\cal C}_{\cal L})$ is labeled by $Q$ in the way explained 
after \kalpar. The group,
\eqn\permu{
\prod_{\alpha=1}^L S_{\ell_{\alpha}},} 
acts naturally on the  
Riemann surfaces by exchanging the $\ell_{\alpha}$ holes that end 
on $\gamma_{\alpha}$. The action of \permu\ lifts to the moduli 
space ${\cal M}_{g,\ell,Q}$ and therefore to the cohomology group 
$H^*({\cal M}_{g,\ell,Q})$. We can then project this cohomology 
group into the subspace which is invariant under the symmetry 
associated to the Young tableaux of $R_1, \cdots, R_L$. This projection 
is made through the operator,
\eqn\bigschur{
{\bf S}_{R_1, \cdots, R_L}=\otimes_{\alpha=1}^L {\bf S}_{R_{\alpha}},} 
where the ${\bf S}_{R_{\alpha}}$ are the usual Schur functors (see 
for example \fh). According to \lmv, the integer invariants $N_{(R_1, \cdots,
R_L),g,Q}$ in \frlinks\  have to be interpreted as  Euler characteristics of
the projected cohomologies, 
\eqn\finaleq{
N_{(R_1, \cdots, R_L),g,Q} =\chi({\bf S}_{R_1, \cdots, R_L} 
(H^*({\cal M}_{g,\ell,Q}))).}   
This gives the geometrical content of the integer coefficients that 
appear in the reformulated polynomial invariants \foot{The above formula 
should hold only up to an overall sign. This is related to the analytic 
continuation that one has to perform in order to compare Chern-Simons 
invariants to enumerative invariants \ov\lm\lmv.}.
   
\medskip\noindent{\bf Remarks}:

\noindent$\bullet$ The equation \multi\ encodes the multicovering and bubbling 
phenomena for the Gromov-Witten invariants of open strings. 
A comparison with the closed string case is very 
illuminating. If $X$ is a Calabi-Yau 
manifold, the Gromov-Witten invariants $N_{\beta}^g$ associated 
to genus $g$ curves in the two-homology class $\beta$ can be 
organized in the Gromov-Witten potential,
\eqn\gvpot{
F(x, t)=\sum_{g \ge 0}\sum_{\beta \in H_2(X)} N_{\beta}^g\, 
x^{2g-2}{\rm e}^{-t \beta}.}
This potential can be rewritten in terms of Gopakumar-Vafa invariants 
$n_{\beta}^g$ \gvm\ as follows:
\eqn\gopavafa{ 
F(x, t)=\sum_{g, \beta} \sum_{d >0} n_{\beta}^g {1 \over d} 
\bigl( 2 \sin (dx/2) \bigr) ^{2g-2} {\rm e}^{-t d \beta}.} 
The Gopakumar-Vafa invariants can be computed in terms of 
Euler characteristics of moduli spaces of holomorphically 
embedded surfaces in the 
Calabi-Yau manifold \gvm\kkv. 
The integer invariants $N_{(R_1, \cdots, R_L),g,Q}$ can then be 
regarded as the open string version of the Gopakumar-Vafa invariants,   
and the relation \multi\ encodes all the multicovering and bubbling 
phenomena associated to the Gromov-Witten invariants for open strings, as 
\gopavafa\ does in the context of closed strings. It is worth 
mentioning that both \multi\ and \gopavafa\ are based on an analysis in terms 
of D-branes, and the Jacobian of the Riemann surface plays a crucial role. 

\noindent$\bullet$ The characterization of the integer invariants given in 
\finaleq\ should be studied more carefully. It was obtained in 
\lmv\ with the simplifying assumption that there is no 
degeneration of the Riemann surface along the moduli. 
Although this assumption gives the right structure of $F(V_1, \cdots,
V_L)$,  the definition of the integer invariants should be analyzed in 
more detail along the lines of \kkv. 

\noindent$\bullet$ As a final remark, it is interesting to 
observe how the structure theorem \multi\ encodes the 
structure of the $1/N$ expansion \largen. If we 
define 
\eqn\relelinks{
f_{(\vec k^{(1)}, \cdots, \vec k^{(L)})}(q, \lambda)
=\sum_{\{ R_{\alpha} \} } 
\prod_{\alpha=1}^L 
\chi_{R_{\alpha}}(C(\vec k^{(\alpha)})) f_{(R_1, \cdots,
R_L)}(q,\lambda)}
and
\eqn\kintslinks{
n_{({\vec k}^{(1)}, \cdots, {\vec k}^{(L)}),g,Q}= 
\sum_{ \{ R_{\alpha} \} } 
 \prod_{\alpha=1}^L
\chi_{R_{\alpha}} (C({\vec k}^{(\alpha)})) 
N_{(R_1, \cdots, R_L),g,Q},} 
one can show \lmv\ that \frlinks\ implies,
\eqn\stflinks{
f_{(\vec k^{(1)}, \cdots, \vec k^{(L)})}(q, \lambda)=\biggl( {
\prod_j (q^{{j\over 2}} -q^{-{j\over 2}})^{\sum_{\alpha=1}^L
k_j^{(\alpha)}}  
\over (q^{{1\over 2}}-q^{-{1\over 2}})^2}\biggr)
\sum_{Q}\sum_{g\ge 0} n_{(\vec k^{(1)}, \cdots, \vec k^{(L)}),g,Q}
(q^{{1\over 2}}-q^{-{1\over 2}})^{2g}\lambda^Q.} 
On the other hand, the ``connected'' invariants defined in \fe\ are
related to  the functions defined in \relelinks\ as follows: 
\eqn\stfk{ 
W^{(c)}_{(\vec k^{(1)}, \cdots, \vec k^{(L)})}(q, \lambda) =
\sum_{d|\vec k^{(\alpha)}} 
d^{\sum_{\alpha=1}^L|\vec k{(\alpha)}|-1} 
f_{(\vec k_{1/d}^{(1)}, \cdots, \vec k_{1/d}^{(L)})}(q^d, \lambda^d).}
Using \stflinks\ and \stfk, it is easy to check that $W^{(c)}_{(\vec k^{(1)}, 
\cdots, \vec k^{(L)})}(q, \lambda)$ has a $1/N$ expansion with the 
structure \largen. From this point of view, the expression in terms 
of the integer invariants $N_{(R_1, \cdots, R_L),g,Q}$ provides 
a resummation of the $1/N$ expansion. 

\newsec{Conclusions and open problems}

In this paper we have described a new set of polynomial invariants for
knots and links  which are closely related to the familiar 
quantum-group polynomial invariants. We have also stated a conjecture 
on their general algebraic structure, and  
described the topological content of their coefficients:  
the integer invariants appearing in the
new polynomials are interpreted as a resummation of the Gromov-Witten
invariants, and are identified in terms of topological properties of the
moduli space of Riemann surfaces with holes embedded in a particular way,
fixed by the knot or link under consideration, into the Calabi-Yau
manifold \rescon\ or resolved conifold.

Up to now, the interpretation of the new invariants in terms of 
enumerative geometry has been fully tested only for the unknot 
\kl\ls. The conjecture \frlinks\ on the structure of the reformulated 
polynomial invariants, however, has been shown to be satisfied in a variety of
cases, and since this structure result is a consequence of the 
geometric formulation, this test can be regarded as a further support 
for this formulation. Unfortunately, not very much is known about the 
properties of quantum-group polynomial invariants for higher
dimensional representations (at least for $SU(N)$ with $N$ generic) 
and no test has been carried out beyond
representations whose Young tableau possesses four boxes. For lower
representations the conjecture \frlinks\ has been tested for a series of 
knots and links \lm\lmv\rama.

Further studies should be done on the topological side to compute the
integer invariants \finaleq. Recent work \kl\ls\ has presented a firm path
towards the computation of these quantities. A good starting point could be the
consideration of torus knots, a case for which a proposal for the
corresponding Lagrangian submanifold is already available \lmv.

Another important issue is the search for more structure. Quantum-group
invariants satisfy skein relations which must have some implications on
the reformulated polynomial invariants. The properties behind these
relations have a different nature than the ones contained in the 
conjecture \frlinks. It would be very important to work out the conditions that
this additional structure imposes on the new integer invariants. In
turn, one should answer also the question about its significance taking
into account their topological origin.

The work summarized in this paper should be extended to take into
consideration quantized universal enveloping algebras different than
$U_q({\rm sl}(N, {\bf C}))$. The extension of \gv\ to other gauge groups 
has been already done in \sv, and it involves non-orientable Riemann 
surfaces in a crucial way. 

Finally, one should rephrase many of the unanswered questions in the
theory of knots and links in terms of these new integer invariants. The
approach certainly opens a new perspective to face these problems.
However, much work is first needed to study these invariants from their
topological origin, obtaining some familiarity with their
properties and developing approaches towards their computation.

\vskip 0.5in 
\vbox{\centerline{\bf Acknowledgments}
\bigskip

We would like to thank Cumrun Vafa for collaboration and useful 
conversations on these topics. We are also grateful to 
Ezra Getzler for explaining to us how to find an 
explicit expression for the reformulated invariants.  
The work by J.M.F.L. was supported in part by DGICYT
under grant PB96-0960 and by Xunta de Galicia under 
grant PGIDT00-PXI-20609. The
work of M.M. is supported by DOE grant DE-FG02-96ER40959.}

\listrefs
\bye